# Simplex Closing Probabilities in Directed Graphs

Florian Unger[*], Jonathan Krebs[†], Michael G. Müller[‡]

May 12, 2022

**Abstract**

Recent work in mathematical neuroscience has calculated the directed graph homology of the directed simplicial complex given by the brain's sparse adjacency graph, the so called connectome. These biological connectomes show an abundance of both high-dimensional directed simplices and Betti-numbers in all viable dimensions – in contrast to Erdős–Rényi-graphs of comparable size and density. An analysis of synthetically trained connectomes reveals similar findings, raising questions about the graphs comparability and the nature of origin of the simplices.

We present a new method capable of delivering insight into the emergence of simplices and thus simplicial abundance. Our approach allows to easily distinguish simplex-rich connectomes of different origin. The method relies on the novel concept of an almost-d-simplex, that is, a simplex missing exactly one edge, and consequently the almost-d-simplex closing probability by dimension. We also describe a fast algorithm to identify almost-d-simplices in a given graph. Applying this method to biological and artificial data allows us to identify a mechanism responsible for simplex emergence, and suggests this mechanism is responsible for the simplex signature of the excitatory subnetwork of a statistical reconstruction of the mouse primary visual cortex. Our highly optimised code for this new method is publicly available.

## 1 Introduction

Topological data analysis (TDA) of undirected simplicial complexes gained from neuronal activity correlation has seen recent success in the field of neuroscience [7, 27]. This approach has limitations, though: undirectedness makes causal investigations difficult and acquiring large scale, yet detailed, in vivo activity data is still infeasible, restricting analysis to in silico simulations.

Analyzing connectomes, i.e., the directed adjacency graph on the level of single neurons and (groups of) synapses, might mitigate both shortcomings. Data is readily available in the form of recent large scale sparse connectome reconstructions [4, 16]. These reconstructions stochastically generate connectomes from a large (but in comparison to biology still limited) set of biological facts (e.g., by using position and orientation of different neuron types and their respective probabilities to form connections), resulting in networks which match these statistics. Even more excitingly, there has been rapid progress in the field of dense reconstructions [17, 24]. Here, the connectome of small volumes of actual brain tissue is determined. Electron microscopes together with advanced data processing methods achieve a resolution allowing to uncover even synaptic vesicles, thus enabling an accurate reconstruction of networks of neurons and synapses. This allows analyzing biological ground truth data, assuming accompanying development of appropriate methods and algorithms.

---

[*]Technische Universität Graz, Austria, florian.unger@igi.tugraz.at
[†]Friedrich-Alexander-Universität Erlangen-Nürnberg, Germany, jonathan.krebs@fau.de
[‡]Technische Universität Graz, Austria, mueller@igi.tugraz.at

**An Abundance of Simplices in Connectomes** Thinking about graphs through the lens of TDA, one immediate problem arises: in contrast to simplicial complexes gained by activity correlation, connectomic graphs themselves have a rather trivial homology. Thus the approach commonly used is to first construct the so called flag complex, i.e., greedily searching for all (ordered) cliques in the graph and their relation to each other, forming a special case of a simplicial complex. Once the flag complex has been generated, the full toolset of TDA is available, e.g., one can calculate the Betti-numbers in different dimensions to get an idea about the topology of the whole structure [14, 22], or perform $Q$-analysis to look for connected pathways [23]. The flag complex also enables our new method introduced here, the simplex completion probability.

The approach considering directed flag complexes of connectomes has been pioneered in [22]. There, two main sources of data were investigated: the dense reconstruction of the neuronal network of the roundworm Caenorhabditis Elegans (C.Elegans) [5, 28] as well as the Blue Brain Project's (BBP) sparse statistical reconstructions of the somatosensory cortex of a rat [16, 21]. In this context, the neurons form vertices and the chemical, directed synapses between two different neurons form edges (more precisely, the usually present multiple synapses from one neuron to another are grouped into a single edge), forming a simple directed graph.

One of the main findings during the investigation of C.Elegans and BBP was an abundance of simplices: there were billions of acyclic, fully connected (i.e., at least one edge between two vertices) subgraphs, way more than in a comparable Erdős–Rényi graph (ER-graph) of equal size and comparable edge density [22].

Inspired by these findings, we investigated connectomes of various forms of networks, such as the Allen Institute's sparse statistical reconstruction of a mouse primary visual cortex [4] and spiking neural networks trained with rewiring [1] and a supervised or unsupervised learning algorithm. All of them similarly showed *far* more simplices than a comparable ER-graph, despite their considerable differences in the nature of origin. Throughout this paper, we will informally refer to this phenomenon as the connectome being simplicially enriched.

**Investigating Simplicial Enrichedness** These findings of strong simplicial enrichedness are intriguing by themselves: they suggest the presence of a strong nonrandom element in the emergence of these graphs. Furthermore, a lack of simplices would render topological analysis of the flag complex rather pointless while simultaneously, a rich presence of high-dimensional simplices motivates a thorough topological analysis. This is enough to ponder about the origin of these simplices, raising the following questions:

1. Can we quantify the idea of simplicial abundance better? Can we refine

the qualitative statement "This graph gives rise to more simplices than a comparable ER-graph" to a quantitative statement for better comparability between simplicially enriched graphs?

2. Can we find an explanation for these high numbers of simplices and quantify the contribution of known principles to simplicial emergence?

3. Can we then efficiently sample new random graphs exhibiting a comparable statistical signature?

In this paper we describe a mathematical and algorithmic framework allowing to tackle these questions. The key idea is the introduction of the so called almost-$d$-simplex, that is, a subgraph with all the structure present to have a $d$-simplex, just short of one edge.

**Why almost-$d$-simplices in particular?**

The particular motif of an almost-$d$-simplex has the following advantages:

- Almost-$d$-simplices and $d$-simplices are closely related: A $d$-simplex is (by construction) just an almost-$d$-simplex together with the missing edge. Due to the combinatoric nature of simplices however an almost-$d$-simplex is not much more than two ($d$-1)-simplices which share a boundary.

- We can measure the ratio of closed or completed almost-$d$-simplices (those almost-$d$-simplices where the missing edge is actually part of the overall graph) versus all almost-$d$-simplices. This measure $p_d \in [0, 1]$ is calculated for every simplex dimension $d \in \{1, \ldots, D\}$, resulting in a vector $p \in [0, 1]^D$. As $p_1$ is nothing more than the well-known edge density of a graph, this could be interpreted as the simplicial generalization of the edge density.

- The combinatorial nature of almost-$d$-simplices allows us (under assumptions) to split contributions to simplex closure by dimension: as Figure 1 illustrates, closing a high-dimensional almost-$d$-simplex simultaneously implies closing lower-dimensional almost-$d'$-simplices (for all $d' \in \{1, \ldots, d-1\}$. In reverse, from a stochastic perspective, higher-dimensional simplex closure can be partially explained by lower dimensional ones. The so called $\hat{p}_d$ formalizes this idea by associating a simplex closure contribution to each dimension $d$, based on the heuristically motivated assumption that "simplex closing is more likely to be caused by low-dimensional rules than high-dimensional ones". This results in ER-graphs having $\hat{p}_d \simeq 0$ for all $d > 1$, making them easily distinguishable from this perspective. Furthermore, mechanisms that could lead to increased 2-simplex emergence like a combination of Hebbian learning and synaptic rewiring can now be quantified as well.

- Almost-$d$-simplices can be enumerated asymptotically optimal, i.e. linear in the number of almost-$d$-simplices. This leads to quick algorithms in practise.

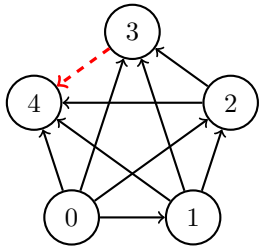

Figure 1: An Almost-4-simplex and the Combinatorics of Closing it. The almost-4-simplex $(0, 1, 2, 3, 4)$ would be closed by the red, dashed edge. Simultaneously, this edge would close one 1-simplex $((3, 4))$, three 2-simplices $((i, 3, 4)$ for $i \in \{0, 1, 2\})$ and three 3-simplices $((0, 1, 3, 4), (0, 2, 3, 4)$ and $(1, 2, 3, 4))$.

**Overview over this Paper**  First we define the graphs, simplices and flag complexes before defining the theory of almost-$d$-simplices as well as their interaction with full $d$-simplices. This allows to define the dimension-wise simplex completion probability $p$. We then describe the combinatorial relation between simplex closings in different dimensions and describe how to measure the contribution *by dimension* with $\hat{p}_d$. Section 2 is devoted to the development of this theory.

Calculating $p$ and $\hat{p}$ for a given graph efficiently requires some algorithmic thought. Thus, in Section 3, we develop a time-optimal algorithm for finding all almost-$d$-simplices in a simplicial complex and describe our `flagser` [14] based C++ implementation `nads`. Note that the availability of python bindings (installable via `pip install pynads`) allows convenient application of our analysis method to other datasets and problems.

In Section 4 we briefly present results for the first two questions, demonstrating the ability of this approach to quantify simplicial enrichedness. We then sketch a potential reason for 2-dimensional simplex closure: common presynaptic neurons increase the connection probability between two neurons due to Hebbian learning. This well known learning rule, together with synaptic rewiring, provides a partial explain for simplex emergence, as a computational experiment demonstrates. Nonetheless explains this mechanism simplicial enrichedness only for certain networks, leaving the main driving force still unclear for others.

In the final section we discuss the benefits and caveats of our approach, sketch potential open questions and research directions, and compare it to other approaches of network analysis.

## 2 Almost-$d$-simplices and their Closing Probabilities

### 2.1 Graphs, Simplices and Flag-Complexes

We follow [22], [23] for the fundamental, preexisting definitions.

**Definition 1** (Simple Directed Graphs)**.** A simple directed graph is a pair $\mathcal{G} = (V, E)$ of a finite set of vertices $V$ and a relation $E \subseteq (V \times V) \backslash \Delta_V$, where

$$\Delta_V = \{(v,v) | v \in V\}.$$

The relation $E$ is the set of directed edges between vertices. The edges are unique ordered pairs $(v,v')$ where we allow reciprocal edges $(v,v')$ and $(v',v)$ to occur simultaneously. However, by excluding $\Delta_V$ from $E$, we exclude self-loops from $v$ to itself. We will write $v \in \mathcal{G}$ or $(v,v') \in \mathcal{G}$ to denote that $v \in V$ respectively $(v,v') \in E$.

As this work concentrates solely on the case of simple directed graphs, we will refer to them simply as graphs.

**Definition 2** (Abstract Directed Simplicial Complex). An abstract directed simplicial complex (or abstract ordered simplicial complex) over a vertex set $V$ is a collection $S$ of non-empty, finite, totally ordered subsets $\sigma \subseteq V$, which is closed under taking non-empty subsets: If $\sigma \in S$ and $\tau$ is a nonempy, order-preserving subset of $\sigma$, then $\tau \in S$. Elements of $S$ of cardinality $d + 1$ are called abstract simplices of dimension $d$ (or abstract $d$-simplex) and the subset of $d$-dimensional abstract simplices is called $S_d$. Furthermore the non-empty subsets of abstract simplices are called faces. We write $\mathrm{Ver}(\sigma)$ to refer to all 0-dimensional faces (the vertices) of $\sigma$ and $\mathrm{Edg}(\sigma)$ to refer to all 1-dimensional faces (the edges).

An abstract $d$-simplex will usually be described by $(d+1)$-tuples $\sigma = (v_0, \dots, v_d)$, where $v_0, \dots, v_d \in V$. The vertex $v_0$ is called the source and $v_d$ the sink of $\sigma$. With $D$ being the maximal abstract simplex dimension of an abstract directed simplicial complex, $S$ is often given as $S = (S_0, \dots, S_D)$.

Note that the total order on an abstract simplex is local to that abstract simplex. There is no requirement of a total order on $V$. Furthermore, it is possible to have two different abstract simplices $\sigma \neq \sigma'$ with $\mathrm{Ver}(\sigma) = \mathrm{Ver}(\sigma')$. In this case, the only difference would be in the orders associated to $\sigma$ and $\sigma'$. See Example 1 for an example of this in the context of directed graphs.

**Definition 3** (Boundaries and Coboundaries). Faces of an abstract $d$-simplex $\sigma$ which are of dimension $d$-1 are called boundaries of $\sigma$. We have the associated $d + 1$ boundary maps $\partial_i : S_d \to S_{d-1}$ for $i \in \{0, \dots, d\}$ by omitting the $i$th vertex:

$$\partial_i(v_0, \dots, v_d) = (v_0, \dots, v_{i-1}, v_{i+1}, \dots, v_d).$$

The dual definition would be a coboundary: the coboundaries of a simplex $\tau \in S_{d-1}$ are all those $d$-simplices which have $\tau$ as a boundary: $\mathrm{coB}(\tau) = \{\sigma \in S_d \,|\, \exists i \in \{0, \dots, d\}$ such that $\partial_i(\sigma) = \tau\}$.

Note that a coboundary can be specified by a tuple $(i,v)$, describing the position and value of the omitted vertex. However, not every tuple $(i,v) \in \{0, \dots, d\} \times V$ is necessarily a coboundary.

To stress the relation of coboundaries to boundaries we introduce the former already here, despite coboundaries not being used till Section 3.

**Definition 4** (Simplices in Directed Graphs). Let $\mathcal{G}$ be a directed simple graph. A $d$-simplex of $\mathcal{G}$ is a nonempty subgraph of $\mathcal{G}$ consisting of $d+1$ vertices together with a total order $\leq$ represented through the edges via $v \leq v'$ iff $(v,v') \in E$.

A $d$-simplex of $\mathcal{G}$ is usually given as $\sigma = (v_0, \dots, v_d)$ indicating that $(v_i, v_j) \in E$ for all $i, j \in \{0, \dots, d\}$ with $i < j$. An alternative definition of a $d$-simplex in a directed graph would be as a subgraph which is both a ($d$+1)-clique and a directed acyclic graph.

We now have two definitions of simplices: abstract simplices in the context of abstract directed simplicial complexes and simplices of directed graphs in the context of directed graphs. However, with $\sigma = (v_0, \ldots, v_d)$ being a simplex in a directed graph $\mathcal{G}$, all induced subgraphs of $\sigma$ created by taking a nonempty subset of nodes and all remaining edges are again simplices in $\mathcal{G}$. That makes the collection of simplices in a directed graph an abstract direct simplicial complex.

A directed flag complex is where both definitions of simplices align, thus we refer to them simply as simplices for the rest of the paper:

**Definition 5** (Directed Flag Complex). Let $\mathcal{G} = (V, E)$ be a simple directed graph. The directed flag complex $\mathcal{F}l(\mathcal{G})$ is an abstract directed simplicial complex over the vertex set $V$ whose $d$-simplices are $d$-simplices in $\mathcal{G}$.

**Example 1.** We investigate the flag complex over the graph $\mathcal{G} = (V, E)$ with $V = \{0, 1, 2, 3, 4, 5\}$ and $E = \{(0, 1), (0, 2), (1, 2), (1, 3), (3, 2), (3, 4), (3, 5), (4, 1)\}$ (see Figure below, left). $\mathcal{F}l(\mathcal{G}) = S_0 \cup S_1 \cup S_2$ where $S_0 = V$ and $S_1 = E$. The graph exhibits furthermore the two 2-simplices $(0, 1, 2)$ and $(1, 3, 2)$ which make up $S_2$. Note that while the vertices $\{1, 3, 4\}$ form a clique, the condition of being ordered is violated, as they are cyclic.

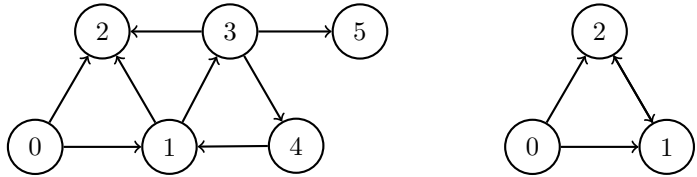

On the right we have another graph, exhibiting two 2-simplices over the same vertex set: $(0, 1, 2), (0, 2, 1)$.

Assuming that every edge on the right hand side goes in both directions, there would be $6 = 3!$ simplices, exactly the number of potential total orders on or permutations of 3 elements. They might not even share a single edge with each other, see e.g. $(0, 1, 2)$ and $(2, 1, 0)$.

This demonstrates how in the directed case of flag complexes, double edges may result in many more simplices over the same number of vertices than the corresponding definitions in the undirected case (or the special case where bidirectional edges are forbidden).

**Remark 1.** Flag complexes are simplicial complexes, but the reverse is not necessarily true. Even if a simplicial complex was defined over a graph with $S_0 = V$ and $S_1 = E$, it does not have to contain every simplex which is a simplex in the graph. In a way, flag complexes can also be understood as the maximal simplicial complex over a given $S_0$ and $S_1$.

## 2.2 Almost-$d$-Simplices

In Section 1 we informally introduced almost-$d$-simplices as simplices missing just a single edge. In case of undirected simplices an almost-$d$-simplex would then just be a pair of ($d$-1)-simplices with a common ($d$-2)-face. In the directed case two complications arise: first, the missing edge direction must be specified as occasionally both are possible. Furthermore, two different simplices over the same set of vertices are now possible (see Example 1), leading to problems as well. Thus we arrive at the following, novel definition:

**Definition 6.** An almost-$d$-simplex is a subgraph of a directed simple graph. For $d \in \{2, 3, \ldots\}$ it is described by a tuple $(\{\sigma, \sigma'\}, e)$ where $\sigma = (v_0, \ldots, v_{d-1})$ and $\sigma' = (v'_0, \ldots, v'_{d-1})$ are $(d-1)$-simplices and $e$ is a tuple of vertices (which may, but does not have to be an edge in the original graph!) such that:

1. There exist $i, i' \in \{0, \ldots, d-1\}$ such that $\partial_i(\sigma) = \partial_{i'}(\sigma')$, i.e. $\sigma$ and $\sigma'$ share a common boundary (in this case, a $(d-2)$-face).

2. $v_i \neq v'_{i'}$, i.e. $\sigma$ and $\sigma'$ differ in (exactly) one vertex.

3. $e$ is either $(v_i, v'_{i'})$ or $(v'_{i'}, v_i)$, where the first is only allowed if $i \leq i'$ and the second only if $i' \leq i$.

We extend the definition to $d = 1$ as $(\{v, v'\}, (v, v'))$ for $v, v' \in V$ with $v \neq v'$.

The extension to $d = 1$ is reasonable as we understand almost-$d$-simplices as a $d$-simplex missing exactly one edge.

A more general definition understanding almost-$d$-simplices not as subgraphs, but objects defined with the means of an abstract directed simplicial complex is certainly feasible. We stick to the more restricted case of graphs and flag complexes as it appropriate for the remaining paper (See Remark 4).

**Example 2.** Here we observe several almost-3-simplices: To the left we have

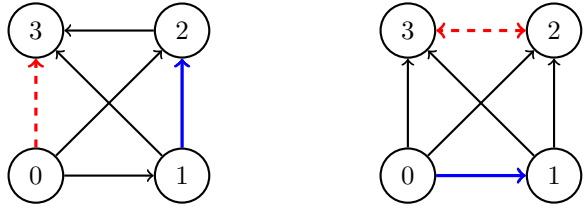

$(\{(0, 1, 2), (1, 2, 3)\}, (0, 3))$. With $i = 0$ and $i' = 2$, the common face is the 1-simplex $(1, 2)$ which is marked blue and thick. The unique "missing" edge in order to form a 3-simplex is red and dashed. Note that its reverse $(3, 0)$ would introduce a cycle 0-1-3-0 into the graph, the reason why it is forbidden. On the right we actually have two almost-3-simplices over the same set of 2-simplices: as $\sigma = (0, 1, 2)$ and $\sigma' = (0, 1, 3)$ over the common face $(0, 1)$ implies $i = i' = 2$, both the edges $(2, 3)$ and its reverse $(3, 2)$ are valid.

**Remark 2.** Let us once again have a look at the the two simplices $\sigma = (0, 1, 2)$ and $\sigma' = (0, 2, 1)$ sketched in Example 1 and check if they form an almost-3-simplex $(\sigma, \sigma', e)$: With $i = 1$ and $i' = 2$ the boundary maps lead to the shared 1-face $\partial_i(\sigma) = (0, 2) = \partial_{i'}(\sigma')$. But $v$ and $v'$ would both be 1, leading to a self-loop $e = (1, 1)$ if not explicitly forbidden by condition ii). Furthermore we simply lack an additional vertex to have proper 3-simplex if $e$ was included.

By design, adding the edge $e$ (which may or may not have been part of the graph already) to an almost-$d$-simplex completes this almost-$d$-simplex to a full $d$-simplex:

**Lemma 1.** *Let $(\{\sigma, \sigma'\}, e)$ be an almost-d-simplex in a graph. Then the flag complex over the graph $\mathcal{G} = (V, E)$ with $V = \mathit{Ver}(\sigma) \cup \mathit{Ver}(\sigma')$ and $E = \mathit{Edg}(\sigma) \cup \mathit{Edg}(\sigma') \cup \{e\}$ contains exactly one d-simplex $\overline{\sigma}$.*

*Proof.* Let the vertices be named as in above definition and let w.l.o.g. $i \leq i'$. Consider
$$\overline{\sigma} = (\overline{v}_0, \ldots, \overline{v}_d) \text{ where } \overline{v}_j = \begin{cases} v_j \text{ for } 0 \leq j < i', \\ v'_{i'} \text{ for } j = i', \\ v_{j-1} \text{ for } i' < j \leq d. \end{cases}$$

We claim that $\overline{\sigma}$ is a $d$-simplex and thus need to check for the existence of edges $\overline{v}_k$ to $\overline{v}_l$ whenever $k < l$. For $k, l \in \{0, \ldots, d\} \setminus \{i'\}$ these are supplied by $\sigma$, as are the edges in case $k, l \in \{0, \ldots, d\} \setminus \{i\}$ supplied by $\sigma'$. The missing edge $(\overline{v}_i, \overline{v}_{i'})$ is exactly $e$. We check uniqueness by showing that there are no more than $\binom{d+1}{2} = \frac{d^2+d}{2}$ edges in $\mathcal{G}$. Adding the edges of $\sigma$ and $\sigma'$, subtracting the shared edges (which form a $(d-2)$-simplex) and adding the single edge $e$ we get:

$$2\frac{(d-1)^2 + (d-1)}{2} - \frac{(d-2)^2 + (d-2)}{2} + 1$$
$$= \frac{(2d^2 - 2d) - (d^2 - 3d + 2) + 2}{2}$$
$$= \frac{d^2 + d}{2}.$$

□

**Definition 7** (Closing or Completing). Let $(\{\sigma, \sigma'\}, e)$ be an almost-$d$-simplex in a graph $\mathcal{G} = (V, E)$. If the "missing" edge $e$ exists (i.e. $e \in E$), the almost-$d$-simplex $(\{\sigma, \sigma'\}, e)$ is completed or closed by $e$ to the $d$-simplex $\overline{\sigma}$ as given by Lemma 1.

Note that this interpretation relies on the greedy flag complex construction, see Remark 4.

The correspondence of almost-$d$-simplices and $d$-simplices is not one-to-one though: One full $d$-simplex contains as many almost-$d$-simplices completing to this particular simplex as it has edges:

**Lemma 2.** *A $d$-simplex $\overline{\sigma}$ contains exactly $\binom{d+1}{2}$ almost-d-simplices which complete to $\overline{\sigma}$.*

*Proof.* Let $\overline{\sigma} = (\overline{v}_0, \ldots, \overline{v}_d)$ and let $0 \leq k < l \leq d$ and thus $e = (\overline{v}_k, \overline{v}_l)$ be one of the $\binom{d+1}{2}$ edges of $\overline{\sigma}$.

We define $\sigma = \partial_l(\overline{\sigma}) = (v_0, \ldots, v_{d-1})$ and analogously $\sigma' = \partial_k(\overline{\sigma}) = (v'_0, \ldots, v'_{d-1})$, where
$$v_i = \begin{cases} \overline{v}_i \text{ for } i < l \\ \overline{v}_{i+1} \text{ else} \end{cases} \quad \text{and} \quad v'_i = \begin{cases} \overline{v}_i \text{ for } i < k \\ \overline{v}_{i+1} \text{ else} \end{cases}$$

It remains to show that with $i = k$ and $i' = l - 1$ the tuple $(\{\sigma, \sigma'\}, (\overline{v}_k, \overline{v}_l))$ is an almost-$d$-simplex:

1. $\partial_i(\sigma) = \partial_i \circ \partial_{i'}(\overline{\sigma}) = \partial_{i'} \circ \partial_i(\overline{\sigma}) = \partial_{i'}(\sigma')$, thus $\sigma$ and $\sigma'$ share a common boundary.

2. By construction $v_i = \bar{v}_i = \bar{v}_k$ and $v'_{i'} = \bar{v}_{i'+1} = \bar{v}_l$. As both are part of the original simplex $\bar{v}$ and $k \neq l$ they must be different.

3. By construction $i \leq i'$, thus $e = (v_i, v'_{i'})$ is permitted.

This construction results in a different almost-$d$-simplex per edge and thus $\binom{d+1}{2} = \frac{d^2+d}{2}$ almost-$d$-simplices per $d$-simplex. Furthermore, possible values of $\sigma, \sigma'$ are completely determined by $e$: vertex-wise, the shared simplex $\check{\sigma}$ of $\sigma, \sigma'$ can only be $\text{Ver}(\bar{\sigma}) \backslash \{\bar{v}_k, \bar{v}_l\}$ and thus the vertex sets of $\sigma$ and $\sigma'$ are given. Due to a lack of additional edges, the vertex order must be inherited from $\bar{\sigma}$ and is thus unique. Therefore more almost-$d$-simplices than edges are impossible. □

## 2.3 The Simplex Completion Probability $p$

An important parameter describing random graphs is the edge density $p_e = \frac{|E|}{|V|^2 - |V|}$. We extend this definition from the simplicial viewpoint to higher dimensions by asking: if we have a graph whose flag complex contains almost-$d$-simplices, what is the chance that the potentially missing edge specified by an almost-$d$-simplex is already part of the graph? Alternatively, what is the change that a random almost-$d$-simplex in the graph is completed? For a concrete graph we can measure it by taking the weighted fraction of completed almost-$d$-simplices with respect to all almost-$d$-simplices.

The last two lemmata justify the combinatoric factor in the following definition:

**Definition 8.** Let $\mathcal{G}$ be a simple directed graph and $S = \mathcal{F}l(\mathcal{G})$ be its directed flag complex. $D$ denotes the maximum simplex dimension, i.e. $D := \max\{n \mid S_n \neq \emptyset\}$. We define the dimension dependent almost-$d$-simplex completion probability for $d \in \{1, \ldots, D\}$ by

$$p_d(\mathcal{G}) := \frac{\binom{d+1}{2}|S_d|}{|\{x \mid x \text{ is an almost-}d\text{-simplex in } S\}|}.$$

The maximum $p_d(\mathcal{G}) = 1$ can be achieved by $\mathcal{G}$ consisting solely of a $d$-dimensional undirected simplex, i.e. vertices $(v_0, \ldots, v_d)$ and edges $(v_i, v_j)$ whenever $i \neq j$. Conversely, the minimum value $p_d(\mathcal{G}) = 0$ can be achieved with $\mathcal{G}$ consisting solely of dijunct almost-$d$-simplices.

**Remark 3.** Note that for a given $D \in \mathbb{N}$ and rational $p^* \in (0, 1)^D$ we may construct a graph $\mathcal{G}$ with $p(\mathcal{G}) = p^*$:

Starting with dimension $D$, we add the desired ratio of disjoint undirected simplices and almost-$d$-simplices (i.e. no overlap in vertices). Note that this does not generate any almost-$(D+1)$-simplices, as they need at least two $D$-simplices which are not disjunct. We then continue one dimension lower, where there are plenty of (almost-)$(D-1)$-simplices from the previous constructions in dimension $D$. Nevertheless, their weighted ratio $p(\mathcal{G})_{D-1}$ is a rational number and thus is the difference $p(\mathcal{G})_{D-1} - p^*_{D-1}$. Thus we apply above construction to the difference. This does not alter $p(\mathcal{G})_D$. Proceed till dimension 1 for a graph with the desired properties.

**Remark 4.** Note that while there is no formal obstruction towards defining $p$ in above fashion on simplicial complexes (instead of flag complexes), the interpretation is not that clear anymore. Not being completed to a full $d$-simplex could have two possible causes: either the missing edge is not part of the simplicial complex, or the simplicial complex simply does not exhibit this particular simplex, even if potential faces would be part of the complex already. The latter aspects are investigated e.g. in [3].

## 2.4 The Simplex Completion Probability Contribution $\hat{p}$

As Figure 1 suggests, the closing of one almost-$d$-simplex has combinatoric effects on lower-dimensional almost-$d$-simplices. This is formalized in the following lemma:

**Lemma 3.** *For all $i \in \{1, ..., d\}$ an almost-$d$-simplex $(\{\sigma, \sigma'\}, e)$ contains $\binom{d-1}{i-1}$ almost-$i$-simplices which can be completed with $e$.*

*Proof.* Let $e := (e_1, e_2)$. Fixing $i$ we first look for the vertex sets of potential shared ($i$-2)-simplices $\breve{\sigma}_i$. As they must be chosen from $(\text{Ver}(\sigma) \cup \text{Ver}(\sigma')) \setminus \{e_1, e_2\}$, we have $\binom{d-1}{i-1}$ choices. Now for $j \in \{1, 2\}$ the order on $\sigma_j$ with $\text{Ver}(\sigma_j) = \text{Ver}(\breve{\sigma}_i) \cup \{e_j\}$ is inherited from the order of $\sigma$ respectively $\sigma'$ and is thus unique. Compliance with the axioms of almost-$i$-simplices for $(\{\sigma_1, \sigma_2\}, (e_1, e_2))$ is clear by construction. □

Another way of stating Lemma 3 would be: "Closing an almost-$d$-simplex simultaneously closes, for $1 \leq i < d$, $\binom{d-1}{i-1}$ almost-$i$-simplices".

As with $p$, we adopt the perspective that our graph was generated by stochastic means and some underlying parameters $p_d$ for $d \in \{1, ..., D\}$. To motivate $\hat{p}$, we go one step further: Simpler rules involving less vertices should be favoured instead of more complicated rules involving many vertices.

The simplest, 1-dimensional rule would then be $p_e$, the overall chance for a potential edge to exist in the graph. We would expect for a simple ER-graph to only exhibit this 1-dimensional contributor. A more complicated rule would involve three vertices (like a combination of Hebbian learning, Dale's law and synaptic rewiring, see Section 4.2.1).

One could now ask: How much does e.g. $p_3$, the 3-simplex closing probability, deviate from the value expected by 1- and 2-dimensional effects? We recursively split contribution to $p$ by the dimension responsible involving the combinatoric factors gained by Lemma 3:

**Definition 9.** We recursively define the $d$-dimensional simplex closing contribution $\hat{p}_d$ for all $d \in \{1, ..., D\}$:

$$\hat{p}_d = p_d - \sum_{i=1}^{d-1} \binom{d-1}{i-1} \hat{p}_i.$$

In particular $\hat{p}_1 = p_1 = p_e$.

Note that $\hat{p} = (\hat{p}_1, ..., \hat{p}_D)$ contains no more additional information about the graph than $p$, it is merely another perspective given by an invertible, linear transformation.

**Remark 5.** Consider an ER-graph generated with edge probability $p_e$ which has an almost-$d$-simplex. The probability $p_d$ that this almost-$d$-simplex is completed to a simplex is $p_e$, as all possible edges are independent.

So $\hat{p}_1 = p_e$ and $\hat{p}_d = 0$ for $d \geq 2$ by induction:

Assume $\hat{p}_i = 0$ for $i \in \{2...(d-1)\}$. In the base case $d = 2$, the assumption is empty and thereby true. Stepping from $d - 1 \mapsto d$ we have:

$$\hat{p}_d = p_d - \binom{d-1}{0}\hat{p}_1 - \sum_{i=2}^{d-1}\binom{d-1}{i-1}\hat{p}_i = p_e - p_e - 0 = 0$$

We observe that $\hat{p}_d$ purely splits off the one dimensional contribution of $p_e$ resulting in zero contribution for all higher dimensions.

While $p_d$ is not measured as defined in Definition 8, the experiment shown in Figure 3 confirms this result.

**Remark 6.** For an effect creating a positive $\hat{p}_2$ see Section 4.2 about Hebbian learning and rewiring.

The values of $\hat{p}_d$ can also be negative: Consider a graph formed by disjoint almost-2-simplices. Here the closing of 2-simplices is suppressed, and we have $p_1 = \frac{2}{6}$ while $p_2 = 0$. This $p_2$ is less than expected in a comparable ER-graph, and thus $\hat{p}_2 = -\frac{2}{6}$ is negative.

In Section 4.1.2 we also use a variation of $\hat{p}$, the so defined

$$\hat{p}_1^{(2)} = \hat{p}_1, \quad \hat{p}_2^{(2)} = \hat{p}_2^{(2)} \quad \text{and} \quad \hat{p}_d^{(2)} = p_d - \hat{p}_1 - (d-1)\hat{p}_2.$$

The resulting measure is equal to $\hat{p}$ up to $d = 3$, but for higher dimensions, we only subtract the influence of $\hat{p}_1$ and $\hat{p}_2$. This representation emphasizes the information "What are the simplex closing probabilities which can not be explained by $p_e$ and a potentially increased almost-2-simplex closing probability?".

## 3 Algorithmic Aspects of Almost-$d$-Simplices

Calculating $p$ has two challenges: Counting, for all $d$, all $d$-simplices and counting all almost-$d$-simplices. For the former, available software such as `flagser_count` [14] solves the problem outstandingly well already. For the latter, we describe an algorithm here.

To keep the following equations brief, we use in this section the shorthands $N_d = |S_d|$ for the number of $d$-simplices and for the number of almost-$d$-simplices we use $N_d^A = |\{a \mid a \text{ is an almost-}d\text{-simplex in } S\}|$.

### 3.1 Enumerating all Almost-$d$-Simplices

Let us assume that we already have a simplicial complex (of a graph), given by arrays $S_0, \ldots, S_D$, where $D$ is the maximum simplex dimension. Vertices are enumerated and mapped to positive integers, simplices are then simply given by a list/tuple of vertices.

Conceptually, to find all almost-$d$-simplices, we start with a shared $(d-2)$-simplex $\breve{\sigma}$. Then we test all pairs of coboundaries $\sigma, \sigma'$ for the necessary conditions stated in Definition 6. By iterating over all $\breve{\sigma} \in S_{d-2}$, we find all almost-$d$-simplices. See Figure 2 for a visual aid.

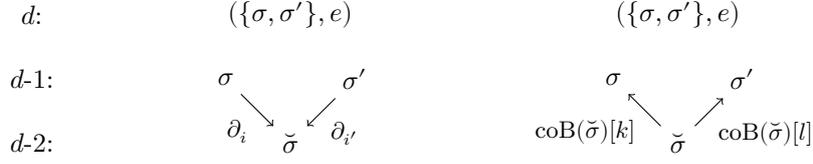

Figure 2: On the left side how one would look for almost-$d$-simplices in the spirit of Definition 6: Starting at $\sigma, \sigma' \in S_{d-1}$ one looks for a common boundary $\breve{\sigma} \in S_{d-2}$ and checks the axioms along the way. On the right hand we see the algorithmic approach utilized in the following code: One starts at $\breve{\sigma} \in S_{d-2}$ and then, utilizing the coboundaries, goes upwards to $\sigma, \sigma' \in S_{d-1}$.

#### 3.1.1 Calculating Cobboundaries and Enumerating almost-$d$-simplices for a given $d$

First we state an time-efficient algorithm for enumerating all the coboundaries from all ($d$-1)-simplices to all $d$-simplices. Note that in our actual implementation, we use a slightly modified algorithm for better space efficiency (see 3.2):

---
**Algorithm 1:** coboundaries($d, S_d, S_{d-1}$)

---
**Output:** A lookup table coB which, for each $\tau \in S_{d-1}$, offers a dynamic array of its coboundaries

1: coB ← new lookup table of size $N_{d-1}$
2: **for** $\tau \in S_{d-1}$ **do**
3:     coB($\tau$) ← new dynamic array
4: **for** $\sigma = (v_0, ..., v_d) \in S_d$ **do**
5:     **for** $i \in \{0, ..., d\}$ **do**
6:         coB($\partial_i(\sigma)$).append(($\sigma, i, v_i$))
7: **return** coB

---

The correctness of this algorithm follows directly from Definition 3. The runtime lies in $\Theta(N_{d-1} + (d+1)N_d)$. We refer to the coboundaries of $\tau$ by coB($\tau$). They are stored in a zero-indexed (dynamic) array accessed on index $k$ via coB($\tau$)[$k$].

For $d \in \{2, ..., D+1\}$ the following algorithm enumerates all almost-$d$-simplices:

**Algorithm 2:** ads_enumerator$(d, S_{d-1}, S_{d-2})$

**Output:** A (dynamic) array of almost-$d$-simplices
1: coB $\leftarrow$ coboundaries$(d-1, S_{d-1}, S_{d-2})$
2: ads $\leftarrow$ new dynamic array
3: **for** $\breve{\sigma} \in S_{d-2}$ **do**
4:    **for** $0 \leq k < l < |\text{coB}(\breve{\sigma})|$ **do**
5:       $(\sigma, i, v) \leftarrow \text{coB}(\breve{\sigma})[k]$
6:       $(\sigma', i', v') \leftarrow \text{coB}(\breve{\sigma})[l]$
7:       **if** $v \neq v'$ **then**
8:          **if** $i \leq i'$ **then**
9:             ads.append$((\{\sigma, \sigma'\}, (v, v')))$
10:         **if** $i \geq i'$ **then**
11:            ads.append$((\{\sigma, \sigma'\}, (v', v)))$
12: **return** ads

**Correctness** To verify correctness we need to ensure that `ads_enumerator` neither misses an almost-$d$-simplex nor counts one multiple times or wrongly identifies something which is not an almost-$d$-simplex as such. Lines `1:`, `5:` and `6:` ensure, together with the correctness of `coboundaries`, that axiom 1 of Definition 6 is adhered to. Line `7:` corresponds directly to axiom 2 and lines `8:` till `11:` correspond directly to axiom 3 in Definition 6, ensuring that we only add almost-$d$-simplices to `ads`.

Furthermore, by line `3:`, we systematically go over all potential shared simplicies $\breve{\sigma} \in S_{d-2}$ and by line `4:` through (up to symmetry) all possible nonequal tuples of $\sigma, \sigma' \in S_{d-1}$.

This way we find all almost-$d$-simplices exactly once.

### 3.1.2 Asymptotic Complexity

**Proposition 1.** *Assuming that a simplicial complex $S$ is already given by arrays $S_0, \ldots, S_D$ and with $d \in \{2, \ldots, D+1\}$ the function call*

$$\text{ads\_enumerator}(d, S_{d-1}, S_{d-2})$$

*enumerates all almost-$d$-simplices in $S$ with an asymptotic time complexity of*

$$\mathcal{O}(N_{d-2} + d^2 N_{d-1} + N_d^A).$$

*Proof.* Line `1:` calls `coboundaries` with a runtime in $\mathcal{O}(N_{d-2} + dN_{d-1})$. Thus we have to prove that the following lines are bounded by $\mathcal{O}(d^2 N_{d-1} + N_d^A)$. As $N_d^A$ is actually the number of almost-$d$-simplices, the runtime of lines `2:` to `11:` are accounted for by that term *as long as the condition on line `7:` is fulfilled*.

It may happen that $\breve{\sigma}$ has two different coboundaries which share the same edge set. This is the exactly when the added vertex is inserted at different spots (see Remark 2) and makes the check at line `7:` fail. We show now that this can only happen up to $d^2 N_{d-1}$ times and is thus bounded by our proposed result as well.

Let $X_d = \{(\sigma, \sigma') \in S_{d-1}^2 | \text{Ver}(\sigma) = \text{Ver}(\sigma'), \sigma \neq \sigma' \text{ and } \sigma, \sigma' \text{ share a common boundary}\}$. The set $X_d$ describes all cases which make it to line `7:`,

but do not fulfill the condition. We now calculate the maximal cardinality of $X_d$: Starting from $\sigma \in S_{d-1}$ there are $d$ choices for a shared boundary $\breve{\sigma}$, as $\breve{\sigma}$ must be $\partial_i(\sigma)$ for some $i \in \{0, \ldots, d-1\}$. Potential coboundaries $\sigma'$ of $\breve{\sigma}$ such that $(\sigma, \sigma') \in X_d$ must yield to the requirement of the same set of vertices of $\sigma$ and $\sigma'$. Thus the to-be-inserted vertex is already set, which leaves only the position to be specified. But there are only $d = |\{0, \ldots, d-1\}|$ choices.

As there are $N_{d-1}$ different choices of $\sigma$, $d$ choices of boundaries from $\sigma$ to $\breve{\sigma}$ and a further $d$ choices from $\breve{\sigma}$ to $\sigma'$ we achieve the desired boundary of $d^2 N_{d-1}$ for $|X_d|$. This shows the desired result. □

The boundary $\mathcal{O}(N_{d-2} + d^2 N_{d-1} + N_d^A)$ is just an upper bound, though. A trivial lower bound would be $\Omega(N_d^A)$, as finding and enumerating all almost-$d$-simplices surely takes as least as long as just enumerating them. Due to the search for coboundaries and wasting time checking the elements of $X_d$, we thus have the asymptotical computational ballast of $\mathcal{O}(N_{d-2} + d^2 N_{d-1})$ which prevents us from directly claiming asymptotic optimality of `ads_enumerator`.

This changes though, once we aim to enumerate all almost-$d$-simplices of *all* dimensions $d \in \{1, \ldots, D+1\}$: As Lemma 2 shows us that $N_d^A \geq \binom{d+1}{2} N_d$ we can hide the unwanted computational ballast $N_{d-2}$ and $d^2 N_{d-1}$ asymptotically behind $N_{d-2}^A$ respectively $N_{d-1}^A$. This is the core idea of the following statement:

**Theorem 1.** *Assuming that a simplicial complex is already given by arrays $S_0, \ldots, S_D$, the algorithm*

---
**Algorithm 3:** `all_ads_enumerator`$(S_0, S_1, \ldots, S_D)$

**Input:** A simplicial complex $S$ in form of arrays $S_0, \ldots, S_D$
**Output:** An array of dynamical arrays containing all
almost-$d$-simplices for $d \in \{1, \ldots, D+1\}$

1: a1s = `new dynamic array`
2: **for** $v, v' \in S_0$ with $v \neq v'$ **do**
3:     a1s.append( $(\{v,v'\},(v,v'))$ )
4: **for** $d = 2, \ldots, D+1$ **do**
5:     a$d$s = `ads_enumerator`$(d, S_{d-2}, \texttt{coboundaries}(d-1, S_{d-1}, S_{d-2}))$
6: **return** a1s, a2s, $\ldots$, a$(D+1)$s

---

*enumerates all almost-$d$-simplices for $d \in \{1, \ldots, D+1\}$ in the asymptotic optimal time*

$$\Theta\left(\sum_{d=1}^{D+1} N_d^A\right).$$

*Proof.* The correctness of lines 2: and 3: follow from the fact that almost-1-simplices are in one-to-one correspondence with edges. Enumerating them costs $\mathcal{O}((N_0)^2) = \mathcal{O}(N_1^A)$ time.

The asymptotic complexity of lines 4: and 5: is

$$\mathcal{O}(\underbrace{N_0 + 1^2 N_1 + N_2^A}_{\text{complexity for } d=2} + \underbrace{N_1 + 2^2 N_2 + N_3^A}_{\text{complexity for } d=3} + \cdots + \underbrace{N_{D-2} D^2 N_D + N_{D+1}^A}_{\text{complexity for } d=D+1}),$$

which together with the cost for the almost-1-simplices boils down to

$$\mathcal{O}(\underbrace{N_0^2 + N_0}_{\in \mathcal{O}(N_1^A)} + \underbrace{\sum_{d=1}^{D}(d^2+1)N_d}_{\in \mathcal{O}(N_d^A)} + \sum_{d=2}^{D+1} N_d^A),$$

where the second estimate is thanks to Lemma 2, yielding our desired result.

The sharpness of the boundary as well as the asymptotic optimality are the same argument: Any algorithm capable of enumerating all almost-$d$-simplices has a lower bound in the number of almost-$d$-simplices. □

**Remark 7.** As checking wether an almost-$d$-simplex completes to a $d$-simplex or not is $\mathcal{O}(1)$ and enumerating almost-$d$-simplices only requires $S_{d-2}$ and $S_{d-1}$, we can efficiently construct the simplicial complex from a directed graph on the fly. This has the same asymptotic bond and is thus optimal as well. Thus the requirement of the flag complex $S$ being computed before starting to enumerate all almost-$d$-simplices could actually be dropped without changing the asymptotic runtime.

## 3.2 `flagser`-based Implementation

Calculating $p_d$ requires just $N_d$ and $N_d^A$ instead of the actual almost-$d$-simplices. So, for our C++ implementation of counting the number of all almost-$d$-simplices `nads` we simply count them instead of enumerating them.

Above algorithms require a flag complex. For our original usecase (analyzing connectomes) `flagser` was the natural choice, as it was engineered exactly for this purpose. It was then natural to hijack `flagsers` data structures and parallisation framework as well and implement our C++ implementation `nads`[1] on top of it. Furthermore, we used `pybind11` [11] in a manner similar to `pyflagser` [6] of giotto-tda [26] to create the `python` bindings `pynads` [2].

Despite just counting almost-$d$-simplices instead of saving them, memory usage is still an issue. So we further adopt and adapt the memory efficient local coboundary search routine of `flagser` to our needs.

---

[1] https://github.com/flomlo/nads
[2] https://github.com/flomlo/pynads – installation e.g. with: `pip install pynads`

|  | Input Size | | Output Size | | Runtime | | |
|---|---|---|---|---|---|---|---|
| **Dataset** | $N_0$ | $N_1$ | $\sum_{d=0}^{D} N_d$ | $\sum_{d=1}^{D+1} N_d^A$ | `flagser_count` | `nads` | `flagser` |
| C.Elegans | 279 | 2194 | $1.5 \times 10^9$ | $1.4 \times 10^{11}$ | 1.3ms | 6.2ms | 57ms |
| BBP | 31K | 7.6M | $1.5 \times 10^9$ | $1.4 \times 10^{11}$ | 7s | 2m35 | 1h10m |
| Mouse V1 | 76K | 7.6M | $4.0 \times 10^8$ | $1.7 \times 10^{10}$ | 11s | 15m45 | >5h |
| q-rewiring | 128 | 4228 | $1.5 \times 10^6$ | $4.7 \times 10^7$ | 12ms | 230ms | 1m47s |

Table 1: Runtimes of `flagser_count`, which generates a simplicial complex from a graph and counts all simplices, `nads`, which counts all almost-$d$-simplices, and `flagser` (with approximation parameter set to 100), which computes the Betti numbers of the simplicial complex up to some approximation. All computations where run on a single 8-core Ryzen 3700 CPU together with 64 GByte of DDR4-3200 RAM.

We calculate the coboundaries of a single ($d$-1)-simplex represented as tuples $(i, v) \in \{0, \cdots, d-1\} \times V$.

With $\text{Out}(v) := \{v' \in V \mid (v, v') \in E\}$ and $\text{In}(v) := \{v' \in V \mid (v', v) \in E\}$:

---
**Algorithm 4:** `local_coboundaries`($\tau \in S_{d-1}$)

---
**Output:** A dynamic array of all coboundaries of $\tau = (v_0, \ldots, v_{d-1})$
1:   $\text{coB}(\tau)$ = new dynamic array
2:   **for** $i \in \{0, \ldots, d\}$ **do**
3:       $C \leftarrow \bigcap_{0 \leq j < i} \text{Out}(v_j) \cap \bigcap_{i \leq j \leq d-1} \text{In}(v_j)$
4:       **for** $v \in C$ **do**
5:          $\text{coB}(\tau).\text{add}(i, v)$
6:   **return** $\text{coB}(\tau)$

---

A vertex $v$ inserted at position $i$ into the simplex $(v_0, \cdots, v_{d-1})$ must have edges $v_j \rightarrow v$ for $j < i$ and $v \rightarrow v_j$ for $j \geq i$. Like in `flagser`, the set of possible extension vertices, $C$ in line `3:`, is calculated using bit operations on rows of the adjacency matrices In/Out.

This is a classic space/time tradeoff. While computing coB for all $\tau \in S_{d-1}$ at once by using the boundary function from $S_d$ is time-optimal, it requires to store all coboundaries of all simplices in $S_{d-1}$ at once. Thus it has $\mathcal{O}(|\text{coB}_{d-1}|) = \mathcal{O}((d+1)N_d)$ space complexity. On the more local implementation we require $\mathcal{O}(N_{d-1}dN_0)$ time in the worst case, but only the absolutely required $\mathcal{O}(\max\{|\text{coB}(\tau)| \mid \tau \in S_{d-1}\})$ in memory.

A further time/space-tradeoff would be not storing the complete simplicial complex in memory, calculating it on demand instead. In our current implementation we choose (using the `IN_MEMORY` compile flag of `flagser`) saving the flag complex as it reduces computational time and has not been a problem in practice so far. But there are no fundamental principles prohibiting this.

Overall the runtime of the analysis is acceptable even on single CPU systems (see Table 1). As enumerating or counting all almost-$d$-simplices can be done in a map-reduce fashion over almost arbitrarily many computers, runtime should not pose an issue. In practise we noticed quite an influence of $D$ on the algorithm runtime.

# 4 Examples of Application

We now present some example analyses of connectomes extracted from biological data and simulated networks using the methods described in the previous sections. The code used to generate the data, perform the analysis and generate the figures is publicly available at `https://gitlab.tugraz.at/6048DFF64EAA81BC/simplex_completion_prob_experiments`.

## 4.1 Connectome Analysis: Quantifying Simplicial Enrichedness

Analyzing connectomes of various origins consistently shows that simplex counts in high dimension are higher than those of matched ER-graphs. We report here on three example connectomes of different size, edge density and nature:

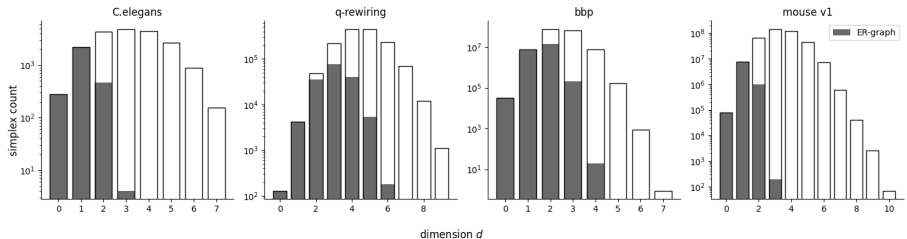

Figure 3: Number of simplices for connectomes of different networks. The simplex counts in dimension *d* are plotted on a log-scale for the network (unfilled) and a comparable ER-graphs (filled).

- *C.Elegans:* the first example is the connectome of the roundworm C.Elegans [5, 28] with 279 neurons and 2194 directed chemical synapses (with a resulting edge density of $p_e \simeq 2.8\%$). This network is completely determined by genetics.

- *BBP:* this connectome comes from the Blue Brain project's sparse statistical reconstruction of the somatosensory microcircuitry of a rat [16, 21] with 31346 neurons and roughly 8 million connections ($p_e \simeq 0.78\%$).

- *Mouse V1:* we further analyzed data from a newer statistical reconstruction of the primary visual cortex of a mouse [4]. Note that we generated a smaller network than the full statistical reconstruction (with $\frac{1}{3}$ of the number of neurons) in order to reduce the computational complexity. This resulted in a network with roughly 76 thousand neurons and 7.6 million edges ($p_e \simeq 0.13\%$).

- *q-rewiring:* we also analyzed the connectome of an artificial spiking neural network, which was randomly initialized and trained with the q-rewiring algorithm [19]. The network was trained in a supervised fashion to perform anomaly detection of a time series of data recorded by stimulating rodent whiskers[20]. It has 128 neurons and 4228 connections ($p_e \simeq 26\%$).

Figure 3 shows the resulting simplex counts. All of the investigated connectomes exhibited more simplices than comparable Erdős–Rényi graph with the same number of vertices and edges, in particular in high dimensions.

### 4.1.1 In-depth Network Analysis with $\hat{p}$

Although all of analyzed connectomes have higher simplex counts than ER-graphs (Figure 3), the results are not comparable in a straightforward manner. It is not directly clear which of the connectomes under investigation are more "simplicially enriched" than others, but such an analysis is important for an understanding of the origins of the high simplex counts and whether artificial mechanism in models match data from real brains. We thus turn to a more in-depth analysis of the simplex counts and the factors underlying them.

An obvious approach to this problem would be to calculate the absolute simplex density, i.e. to calculate the fraction of simplices found between all possible directed simplices: $\frac{N_d}{\binom{N_0}{d+1}(d+1)!}$. We found this to be not very enlightening:

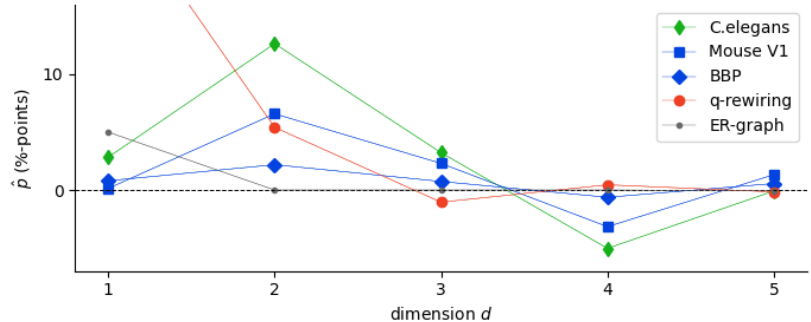

Figure 4: Comparison of simplex closing contributions for different networks. The plot shows the value of $\hat{p}$ up to dimension 5 for the networks shown in Figure 3 an for an ER-graph. Colors indicate the type of connectome (green: biological data, blue: statistical reconstruction, red: artificial network).

When comparing graphs of different size and general edge density, investigating just the raw number of simplices holds no insight on mechanisms leading to more or less simplices compared to pure chance. Furthermore, even recognizing ER-graphs as such is surprisingly difficult – at least compared to an approach utilizing $\hat{p}$, where the signature of ER-graphs is very distinguishable (see Remark 5).

Analyzing the simplex closing contributions $\hat{p}$ (see Section 2) allowed to clearly distinguish the different networks and shed light on the differences between their simplicial structure. The result values of $\hat{p}$ for the networks analyzed above are shown in Figure 4. For a random ER-graph (in this case with $\hat{p}_1 = p_1 = 5\%$ edge density), we have $\hat{p}_d \simeq 0$ for all $d \geq 2$ (as predicted in 5). We find that values of $\hat{p}$ fall between two extremes. First, the artificial, q-rewiring trained network exhibits (beside its outstanding edge density of 26%) a significant almost-2-simplex closing probability. There seem to be no higher-dimensional positive contributions to simplex closing probabilities, though, and a small negative aptitude for 3-dimensional simplices. The other extreme lies in the biological connectome of C.Elegans. While $\hat{p}_2$ dominates, even $\hat{p}_3$ is significantly elevated, highlighting an almost-3-simplex closure aptitude which can not be explained by combinatorial effects of dimensions 1 and 2. The two different statistical reconstructions of rodent brain regions fall between these extremes, but considerably differ from the biological behavior despite their level of biological detailedness.

It is as of yet unclear whether these noticeable differences in the signature of simplex closing contributions are due to differences in reconstruction methods, due to location of the brain regions modeled, or due to other reasons. Ideally, potential mechanisms underlying the differences should be isolated and analyzed, e.g. through simulations. Below, we present a first step in this direction (Section 4.2). Additionally, as data from dense reconstruction (e.g. from mammal cortex) becomes available in reasonable volumes, these should be analyzed and compared to the data from the (arguably rather simple) brain of C.Elegans.

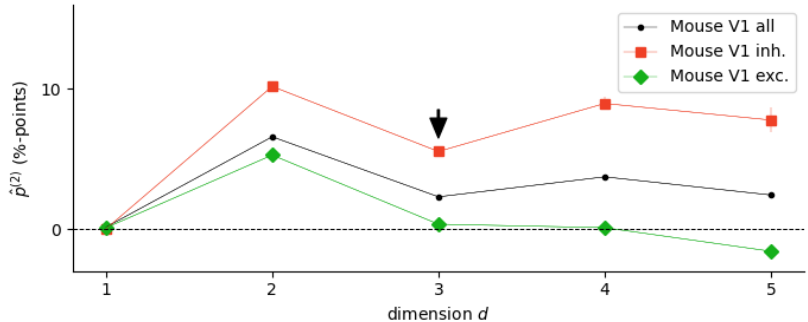

Figure 5: Analysis of $\hat{p}^{(2)}$ for subpopulations of the Mouse V1 connectome. A comparison of $\hat{p}^{(2)}$, the simplex closing contribution adjusted only of 1- and 2-dimensional influences, shows a clear difference between the excitatory (green) and inhibitory (red) subpopulations (black: entire population). Error bars denote standard deviation over $N = 20$ runs (note that in most cases these are too small to be visible).

#### 4.1.2 Analysis of Mouse V1 subpopulations

Biological neurons can be separated into two classes according to their kind of influence on downstream neurons. Excitatory neurons increase, upon spiking, the firing probability of neurons they project to, while inhibitory neurons decrease the probability of firing. As different classes have different properties and connectivity patterns in the brain, they are routine analyzed and modeled separately.

We thus performed an analysis of the simplex closing probabilities for the statistical reconstruction data from mouse V1 (see above) to investigate whether they are similar or different for the excitatory and inhibitory subpopulations. While analyzing both populations together shows $\hat{p}_3 \gg 0$ and thus apparently simplex closing mechanism beyond dimension 2 (see above, Figure 4), the results for individual subpopulations gives a dramatically different result (Figure 5). We here plotted $\hat{p}^{(2)}$ instead of $\hat{p}$ to highlight the differences, i.e., Figure 5 shows (for $d > 2$) the simplex closing probabilities which cannot be explained by an increase in $\hat{p}$ in dimensions 1 and 2 (see above, Section 2.4).

We found that the two subpopulation have radically different $\hat{p}^{(2)}$ signatures. The inhibitory population has already a strong influence of $\hat{p}_2$, but seems to further exhibit a (currently unknown) mechanism which seems to need at least 4 neurons to come into play. The excitatory population behaves radically different: as only $\hat{p}_1^{(2)}$ and $\hat{p}_2^{(2)}$ are significantly above zero, 1- and 2-dimensional simplex closing mechanisms seem to be a suitable explanation for everything.

This naturally raises the question: what might be a suitable 2-dimensional simplex closing mechanism? A candidate is discussed in the next section.

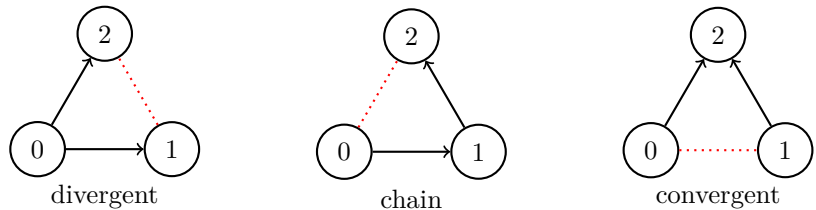

Figure 6: The three constellations (motifs) of almost-2-simplices. The naming follows [18].

## 4.2 Quantifying the Influence of Hebbian Learning and Rewiring on Simplex Emergence

Here next investigate a candidate mechanism underlying simplicial enrichedness of biological connectomes: a combination of Hebbian learning and synaptic rewiring.

### 4.2.1 Hypothesis: Hebbian Learning and Rewiring explain Simplicial Enrichedness

Network activity and structure in biological neural networks is shaped by many different mechanisms. Three important contributors are:

- Hebbian learning [9]: while the biophysical processes shaping synaptic plasticity are complex, this simple mechanism has been found to capture an essential aspect. Hebbian learning is often summarized as "What fires together, wires together", i.e., prolonged correlated firing leads to the strengthening of (existing) synaptic connections. In models, this mechanism is usually amended by some depressing mechanism to prevent unbound growth of all weights.

- Dale's law [25]: neurons have been found to perform the same chemical action at all of their outgoing synapses, i.e., they are either excitatory or inhibitory (see above). In a network model, this corresponds to ensuring that all outgoing weights have the same sign.

- Structural plasticity [8, 10]: synaptic formation and depletion are ongoing processes in the brain, and synaptic lifetime is usually somewhat limited. One kind of modeling of such ongoing rewiring [1, 12] proposes that the overall amount of synapses stays about the same (as they cost energy to maintain). In this method, weak synapses are assumed to be unimportant and are thus occasionally removed, while new potential synapses are randomly drawn and initialized with a small weight.

We propose that all these rules in conjunction would result in an increase in $\hat{p}_2$ for the following reasons. Imagine an excitatory or inhibitory neuron 0 presynaptic to two neurons 1 and 2 (Figure 6, left). By Dale's law 0 has the same kind of influence on 1 and 2, leading to more correlated firing. Now should a random edge sampling process enable a connection $(1,2)$ or $(2,1)$, then this connection would grow stronger by Hebbian learning over time, leading to a stable connection. Thus the almost-$d$-simplices $(\{(0,1),(0,2)\}, e)$ for $e \in \{(1,2),(2,1)\}$ have an increased chance of being closed by this purely local phenomenon. For a second constellation (Figure 6, middle) consider 0 being presynaptic to 1 which is already presynaptic to 2, but no connection between 0 and 2, yet. As long as 0 and 1 are excitatory, there is again a higher correlation in firing and analogous to the argument before an increased chance too complete this almost-2-simplex. Only the right case would not be closed with an increased likelihood by above arguments. In addition to this increase of 2-simplices due to Hebbian learning and Dale's law, ongoing network rewiring in conjunction with Hebbian learning will result in the depletion of connections in which neuron's fire asynchronously, while synaptic connections between neurons exhibiting correlated firing will be stabilized over time (once they have been formed randomly).

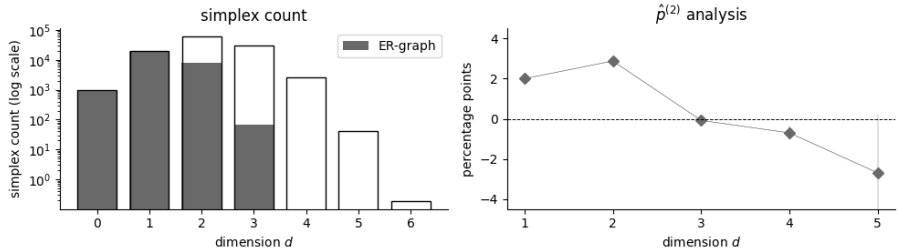

Figure 7: Analysis of spiking neural networks trained with Hebbian learning and rewiring. Left: simplex count (mean values for $N = 100$ trained networks) compared to equivalent ER-graphs. Note that for $d \in \{0, 1\}$, the values of network and ER-graph match. Right: $\hat{p}^{(2)}$ analysis. Errorbars denote standard deviation over $N = 100$ runs (note that in most cases these are too small to be visible).

### 4.2.2 Training Spiking Neural Networks with Rewiring and Hebbian Learning leads to increased $\hat{p}_2$

To test the hypothesis proposed in the previous section we simulated spiking neural networks incorporating these three mechanisms. The networks were run for some period of time, allowing synapses to form and stabilize through the combination of Hebbian learning and rewiring. This allows us to see if this combination of plasticity mechanisms, through shaping the connectome, gives rise to more simplices. Furthermore, by investigating whether $\hat{p}^{(2)} \gg 0$ for these networks, we can test whether this setup is enough to explain all heightened simplex closing probabilities in higher dimensions.

We initialized the connectome of these networks with a self-loop-free ER-graph of 1000 vertices and 50000 potential edges, corresponding to 1000 excitatory neurons. Synaptic weights associated to the edges were always positive (following Dale's law) and were randomly initialized with small values (see the supplementary information for details). Inhibition was modeled in an abstract fashion to keep the network to keep the network activity in a reasonable regime. We ran the simulation until a certain fraction of weights ($p_1 = 0.02$) stabilized, i.e., the weights reached high values due to Hebbian learning. We then dropped weights with low values and analyzed the resulting connectome.

Figure 7 shows that these networks exhibit substantially more simplices than a comparable ER-graph. We furthermore analyzed $\hat{p}^{(2)}$ and found that there seem to be only 1- and 2-dimensional simplex closing effects. In higher dimensions, we surprisingly found negative closing probability contributions, although the variance increased with the dimension.

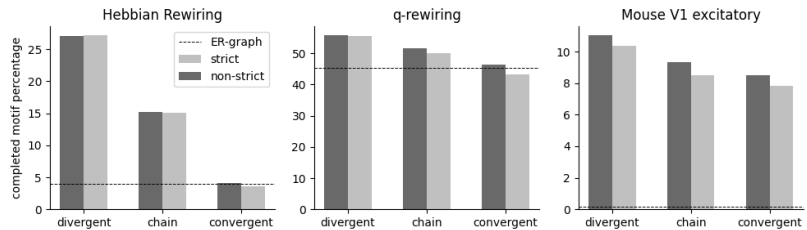

Figure 8: Analysis of motif closing ratios for the motifs shown in Figure 6. The investigated networks all have $\hat{p}_d \simeq 0$ for $d \geq 3$ (see Figures 4, 5, and 7). Left: spiking neural network trained with Hebbian learning and rewiring. Middle: spiking network trained with q-rewiring. Right: excitatory subnetwork from mouse V1.

This shows that Hebbian learning (together with network rewiring) indeed increases $\hat{p}$ for $d = 2$. Our results are qualitatively similar to the findings of the excitatory subpopulation of mouse V1 (large $\hat{p}^{(2)}$ for $d = 2$, but small or slightly negative values for $d \geq 3$). While the reasons for the negative contributions in higher dimensions is are as of yet unclear, our findings suggests Hebbian learning as a prominent factor underlying the structure of excitatory networks in the brain. It has been suggested that Hebbian learning also plays an important role in shaping inhibitory connectivity (e.g., [29]). The difference of the simplex closing signature within inhibitory networks (see above, Figure 5) suggest that other, to date unknown, factors must be considered.

### 4.2.3 A Closer Look: Hebbian Learning and Rewiring Motif Analysis

We have shown that a combination of Hebbian learning, Dale's law and rewiring mechanisms leads to increased $\hat{p}_2$ and keeps higher-dimensional contribution equal or even lower than zero, resulting in a simplex signature matching several examined connectomes including the excitatory subnetwork of Mouse V1. To elucidate the emergence of this simplex signature, we now turn again to the three simple motifs shown in Figure 6. Should our hypothesis on Hebbian learning and rewiring being the main driver hold true, the divergent and chain motifs should exhibit a higher fraction of completed motifs than pure chance level, while the convergent motif should not be affected by this (see above).

To investigate whether this is indeed the case, we determined the ratios $\frac{\#\text{completed } m}{\#m}$ and $\frac{\#\text{completed strict } m}{\text{strict } m}$ for $m \in \{\text{divergent}, \text{chain}, \text{convergent}\}$ (see Figure 6). For this analysis, a motif is considered completed if one or both edges between the neurons connected by the red dashed line exists. We distinguished between a strict analysis, in which motifs are only taken into account if no reciprocal edges exists, and a non-strict analysis, which all motifs are counted. These closing ratios can be compared to the chance level, i.e. the probability of these motifs to be closed in an ER-graph with edge density $p_e$, which is $2p_e - p_e^2$.

The results of this analysis are shown in Figure 8 for three different kinds of network. For the spiking network with Hebbian learning and rewiring, the closing ratios match our theoretical prediction. Interestingly, the network trained with q-rewiring shows similar results. This surprises, as this network was not trained in with Hebbian learning, but in a task-dependent supervised manner using error backpropagation. The closing ratios for mouse V1 differ from the theoretical prediction (based on Hebbian learning). While the divergent motif is still the most likely to be closed, all closing probabilities much higher compared to chance level than expected. Possible explanations for this include local clustering, which would alter the chance level (although in this case the exact value is hard to determine), or that the layer-wise structure of this network together with many lateral connections in the same layer lead to completely different simplex emergence mechanisms.

# 5 Summary and Outlook

Using the concepts of almost-$d$-simplices and simplex closing probabilities, we have developed the theory of a new method to analyze directed graphs. With the accompanying implementation we deliver not only a method, but a practical tool. In the case of graphs whose local structure is assumed to be dominated by local and stochastic rules (e.g., connectomes) this tool helps in achieving a deeper understanding of these local rules, and raises new questions and promising research directions.

Nevertheless, we cannot claim to have fully answered the first two questions raised in Section 1: there is no absolute, strict relation between $p_d$ and the number of $d$-simplices, as ($d$-1)-simplices can arrange in different ways, forming more or less almost-$d$-simplices. A more in-depth analysis of the emergence of almost-$d$-simplex could clarify this relation, but is still an open problem. However, this caveat does not pose a significant problem for our application.

While the inverse problem, i.e., finding graphs which exhibit a given $p$-signature has been solved in principle by Remark 3, sampling graphs where furthermore $N_0$ is fixed can be quite hard. Graphs such as this, i.e., simplicially enriched in a certain way w.r.t. $p$, would serve at least two purposes. First they can serve as null-models to investigate the statistical significance of findings on the simplicial complex, e.g., Betti-numbers or $q$-pathways. Furthermore, they could prove to be useful as initializers for network models such as Liquid State Machines [15], allowing to test the relation between simplicial enrichedness and computational capabilities. Preliminary results regarding the construction of simplicially enriched graphs (not shown) suggest that is may be possible by searching and closing almost-$d$-simplices, but we leave this problem to future work.

We have focused here on directed graphs. However, the generalization to undirected graphs and flag complexes is straightforward. The definition boils down to requiring a pair of nonequal ($d$-1)-simplices sharing a common ($d$-2)-simplex and all the results following can be adapted without much effort.

Another popular network analysis approach targets motifs [2, 3, 13], i.e., induced subnetworks with three, four or even more vertices and a particular structure. Naturally, this approach reveals more details than simplex analysis, but is problematic in its scaling with dimensionality: for motifs consisting of only five vertices (and their possible single-/double-edge configurations), there are already thousands of possible structures.

Simplices are a special subgroup of motifs with two advantages: There is only one motif per dimension and they can be efficiently searched for, which could otherwise pose a problem in large-scale graphs. Additionally there is the entire field of TDA eager to uncover additional insights hiding in the flag complex. Insights gained in that way usually span many high dimensional simplices and would not be discovered by general motif-methods alone.

Almost-$d$-simplices are slightly more intricate than simplices: there are $\binom{d+1}{2}$ different almost-$d$-simplices per dimension if we understand them as a family of motifs. From this perspective, the level of detail of simplex completion probability analysis lies somewhere between counting all simplices and motif-based approaches. From the computational cost it settles neatly between counting simplices and more traditional flag complex analysis tools (see Table 1). We thus believe our method can serve an important role in elucidating the structure of large graphs, and (due to the publicly available implementation) propose it as a low-cost analysis tool.

**Acknowledgements**

We would like to thank Robert Legenstein for his general council, Horst Petschenig for supplying us with q-rewiring trained connectomes, Michael Kerber for the discussion leading to Theorem 1 and two anonymous reviewers for their sizeable contributions to mature this paper.

# Simplex Closing Probabilities in Directed Graphs

## *Supplementary information*

Florian Unger, Jonathan Krebs, Michael G. Müller

April 29, 2022

## Details to Hebbian learning simulations (Section 4.2.2)

**Neuron model**

The network consisted of $N = 1000$ leaky integrate-and-fire (LIF) neurons. Each neuron's membrane potential was updated according to

$$\tau_\mathrm{m} \frac{du}{dt} = -(u - u_\mathrm{rest}) + R_\mathrm{m} (I_\mathrm{syn} + I_\mathrm{inh}) \tag{1}$$

where $\tau_\mathrm{m} = R_\mathrm{m} C_\mathrm{m}$ is the membrane time constant (product of the membrane resistance $R_\mathrm{m}$ and the membrane capacitance $C_\mathrm{m}$), $u_\mathrm{rest} = -65\,\mathrm{mV}$ is the resting potential, $I_\mathrm{syn}$ is the synaptic input current from other neurons, and $I_\mathrm{inh}$ is the inhibitory input current. The values of $\tau_\mathrm{m}$ were randomly drawn for each neuron from a Gaussian distribution with mean $15\,\mathrm{ms}$ and standard deviation $5\,\mathrm{ms}$, clipped to the range $[5\,\mathrm{ms}, 25\,\mathrm{ms}]$. The values of $C_\mathrm{m}$ were similarly drawn from a Gaussian distribution with mean $250\,\mathrm{pF}$ and standard deviation $50\,\mathrm{pF}$, clipped to the range $[50\,\mathrm{pF}, 450\,\mathrm{pF}]$. At time $t_0$ of the discrete-time simulation (using a time resolution of $\Delta t$), a nonrefractory neuron created a spike with probability $\rho(t_0)\Delta t$, with

$$\rho(t) = \frac{1}{\Delta t} \exp\left(\frac{u(t) - u_\mathrm{T}}{T}\right), \tag{2}$$

where we used $u_\mathrm{T} = -55\,\mathrm{mV}$ and $T = 1\,\mathrm{mV}$. If a spike was generated, neurons' membrane potentials were reset to $u_\mathrm{reset} = -75\,\mathrm{mV}$ and they entered a refractory prediod (duration $3\,\mathrm{ms}$) during which $u(t)$ remained clamped to this value.

The inhibition was modeled in an abstract fashion and chosen to keep each neuron's firing rate in a reasonable range (similar to a homeostatic mechanism, but acting on a faster timescale). Each neuron's current rate was estimated by low-pass filtering its spikes (time constant: $500\,\mathrm{ms}$), resulting in an instantaneous rate estimate $r(t)$ for each neuron. The inhibitory current was set to $I_\mathrm{inh}(t) = -a_\mathrm{inh} (r(t) - r_\mathrm{target})$ with $a_\mathrm{inh} = 25\,\mathrm{pA}$ and $r_\mathrm{target} = 10\,\mathrm{Hz}$.

The input current $I_\mathrm{syn}$ of each neuron was generated from incoming spikes via

$$\tau_\mathrm{syn} \frac{dI_\mathrm{syn}}{dt} = -I_\mathrm{syn} + \sum_{j \in \mathrm{PRE}} w_j\, \delta\left(t - t_j^{(f)}\right), \tag{3}$$

where $\tau_\mathrm{syn} = 3\,\mathrm{ms}$ is the synaptic time constant, and the sum runs over the set of presynaptic neurons $j$ with spike times $t_j^{(f)}\ \forall f$ and synaptic weight $w_j$.

**Synapse model**

Synaptic weights were sparsely initialized (see below) by drawing them from a lognormal distribution with mean and standard deviation of the underlying normal distribution set to $0\,\mathrm{pA}$ and $1\,\mathrm{pA}$, respectively, and then clipping to the range $[0\,\mathrm{pA}, w_\mathrm{max}]$ with $w_\mathrm{max} = 20\,\mathrm{pA}$. During the simulation, the weights were updated by a Hebbian rule, with each pairing of pre- and postsynaptic spikes with timing difference $\Delta t$ leading to a weight change of $w \leftarrow w + \Delta w(\Delta t)$ with

$$\Delta w(\Delta t) = \eta \exp\left(-|\Delta t|/\tau_\mathrm{Hebb}\right) \tag{4}$$

where we used $\eta = 1\,\mathrm{pA}$ and $\tau_\mathrm{Hebb} = 20\,\mathrm{ms}$. Additionally, each presynaptic spike led to a weight change of $w \leftarrow w - \eta A_-$ with $A_- = 0.35$. After each update, the weights were clipped to the range $[0\,\mathrm{pA}, w_\mathrm{max}]$.

Synaptic delays randomly drawn from the range $[1, 3]\,\mathrm{ms}$. Such a value was drawn each time a connection was generated (even if this connection had existed earlier with a different delay).

## Rewiring

During the simulation, the connectivity was allowed to change using a form of synaptic rewiring. At any point during the simulation, the number of synapses in existence was chosen such that the connectivity was $C_{\text{init}} = 0.05$ (i.e. for $N = 1000$, the number of synapses was held constant at $N^2 C_{\text{init}} = 50000$). The initial connectivity was chosen randomly out of all possible connections (no autapses were allowed at any point). During the simulation, we updated the connectivity every $1\,s$ by grouping the existing synapses into three sets by applying the thresholds $\theta_{\text{drop}} = 0.1$ and $\theta_{\text{stable}} = 0.5$ to the the current weights (altered by the Hebbian learning rule):

- Low weights: connections with $w/w_{\text{max}} < \theta_{\text{drop}}$ were dropped. For each removed synapse, we introduced new weight by randomly choosing a non-existing connection and setting its weight to $\theta_{\text{new}} w_{\text{max}}$ with $\theta_{\text{new}} = 0.1$.

- Intermediate weights: connections with $\theta_{\text{drop}} \leq w/w_{\text{max}} < \theta_{\text{stable}}$ were kept.

- High weights: connections with $w/w_{\text{max}} \geq \theta_{\text{stable}}$ were also kept. We considered these connections to be stable.

We simulated the network until the fraction of stable connections was $C_{\text{stable}} \geq 0.02$. After the simulation, we removed all but the strongest synapses such that the connectivity was exactly equal to $C_{\text{stable}} = 0.02$. This pruned connectome was used for further analysis.

## Simulation details

The model was implemented using Brian2 [1], and simulated using a discrete time step of $\Delta t = 0.1\,\text{ms}$.

## References

[1] Marcel Stimberg, Romain Brette, and Dan FM Goodman. "Brian 2, an intuitive and efficient neural simulator". In: *Elife* 8 (2019), e47314.